\documentclass{article}
\usepackage{amssymb,latexsym}
\usepackage{graphicx}
\usepackage{amsmath}
\usepackage{amsthm}
\usepackage{amscd}
\usepackage{bm}
\usepackage{indentfirst}
\begin{document}
\newcommand{\D}{\displaystyle}
\newcommand{\DF}[2]{\frac{\D#1}{\D#2}}
\def\Xint#1{\mathchoice
{\XXint\displaystyle\textstyle{#1}}%
{\XXint\textstyle\scriptstyle{#1}}%
{\XXint\scriptstyle\scriptscriptstyle{#1}}%
{\XXint\scriptscriptstyle\scriptscriptstyle{#1}}%
\!\int}
\def\XXint#1#2#3{{\setbox0=\hbox{$#1{#2#3}{\int}$}
\vcenter{\hbox{$#2#3$}}\kern-.5\wd0}}
\def\ddashint{\Xint=}
\def\dashint{\Xint-}
\pagestyle{plain}
\newtheorem{theorem}{Theorem}[section]
\newtheorem{definition}{Definition}[section]
\newtheorem{remark}{Remark}[section]
\newtheorem{corollary}{Corollary}[section]
\newtheorem{proposition}{Proposition}[section]
\newtheorem{lemma}{Lemma}[section]
\normalsize
\title{Analysis on Metric Space $\mathbb{Q}$
\thanks{Department of mathematics, Rice University, Houston, TX 77005, U.S.A; weizhu@math.rice.edu}}
\author{Wei Zhu}
\maketitle
\begin{abstract}
In this paper, we show that the metric space
$(\mathbb{Q},\mathcal{G})$ is a positively-curved space (PC-space)
in the sense of Alexandrov. We also discuss some issues like
metric tangent cone and exponential map of
$(\mathbb{Q},\mathcal{G})$. Then we give a decomposition of this
metric space according to the signature of points in $\mathbb{Q}$.
Some properties of this decomposition are shown. The second part
of this paper is devoted to some basic analysis on the space
$(\mathbb{Q},\mathcal{G})$, like the tensor sum and $L^p$ space,
which can be of independent interest. In the end, we give another
definition of derivative for multiple-valued functions, which is
equivalent to the one used by Almgren. An interesting theorem
about regular selection of multiple-valued functions which
preserves the differentiability concludes this paper.
\end{abstract}
\tableofcontents
\section{Introduction}
The theory of multiple-valued functions in the sense of Almgren
$\cite{af}$ has several applications in the framework of geometric
measure theory. In deed, multiple-valued functions give a very
useful tool to approximate some abstract objects arising from
geometric measure theory. For example, Almgren (see $\cite{af}$)
used multiple-valued functions to approximate some rectifiable
currents, hence successfully got the partial interior regularity
of area-minimizing integral currents. Solomon (see $\cite{sb}$)
succeeded in giving proofs of the closure theorem without using
the structure theorem. His proofs rely on various facts about
multiple-valued functions. There are also some other work
concerning multiple-valued functions, see
$\cite{dgt},\cite{gj},\cite{lc1},\cite{lc2},\cite{mp},\cite{zw1},\cite{zw2}$.
All these work raises the need of further studying of
multiple-valued functions.\\
In $\cite{af}$, Almgren gave an explicit bi-Lipschitzian
correspondence between the metric space $(\mathbb{Q},\mathcal{G})$
and a finite polyhedral cone $\mathbb{Q}^*$ in some higher
dimensional Euclidean space. His analysis on multiple-valued
functions are mainly based on this correspondence. In this paper,
we are focusing on the metric space $(\mathbb{Q},\mathcal{G})$
itself. By carefully studying the geodesics connecting any two
points in $\mathbb{Q}$, we are able to claim:
\begin{theorem}
The space $\mathbb{Q}$ with metric $\mathcal{G}$ is
positively-curved (a PC-space) in the sense of Alexandrov.
\end{theorem}
We also give an explicit description of the abstract tangent cone
in this metric space $(\mathbb{Q},\mathcal{G})$:
\begin{theorem} For any point $A\in\mathbb{Q}$, such that
$S(A)=(J,k_1,\cdot\cdot\cdot,k_J)$,
$$\bold{Tan}_A(\mathbb{Q})=\mathbb{Q}_{k_1}(\mathbb{R}^n)\times
\mathbb{Q}_{k_2}(\mathbb{R}^n)\times \cdot\cdot\cdot\times
\mathbb{Q}_{k_J}(\mathbb{R}^n)$$ with the product metric.
\end{theorem}
The second part of this paper is based on the following definition
of tensor sum of multiple-valued functions. (Recall that there is
no suitable notion of ``addition'' for arbitrary two
multiple-valued functions)
\begin{definition}
Suppose $f(x)=\sum_{i=1}^p [[f_i(x)]]$, $g(x)=\sum_{j=1}^q
[[g_j(x)]]$, where $p$ and $q$ are not necessarily the same.
Define
$$(f\oplus g)(x)=\sum_{i,j} [[f_i(x)+g_j(x)]].$$
(i.e the tensor sum is a $pq-$valued function). \end{definition}
This definition is of limited uses in the sense that even if
$p=q=Q$, the tensor sum of two $Q-$valued functions gives a
$2Q-$valued functions, which makes it hard to talk about
derivatives. We are expecting some good notion of ``addition''
which will enable us to define various things like integration,
differential equation in the setting of multiple-valued
functions.\\
The last part of this paper gives another definition of
derivatives for multiple-valued functions if a priori the function
is continuous. Unlike using linear approximation to define
derivative in $\cite{af}$, which avoids ``subtraction'', our
definition is more calculus-oriented:
\begin{definition} Suppose
$f:\mathbb{R}^m\to\mathbb{Q}$ is continuous, fix
$x_0\in\mathbb{R}^m$ and $v\in\mathbb{R}^m$, the directional
derivative of $f$ at $x_0$ in the direction $v$ is the following
limit if it exists:
$$L(v):=\lim_{t\to 0}\DF{f(x_0+tv)(-)f(x_0)}{t}.$$
\end{definition}
We will show that this definition is equivalent to the one in
$\cite{af}$.\\
We conclude this paper by a selection theorem for multiple-valued
functions. One of the motivations for this kind of question is
whether we can decompose a $Q-$valued function into $Q$
single-valued functions which preserve some properties of the
original function. There are already some important work by De
Lellis, Grisanti and Tilli $\cite{dgt}$ and by Goblet
$\cite{gj}$.\\ As far as continuity is concerned, we have some
positive and also some negative results, namely, a continuous
$\mathbb{Q}_Q(\mathbb{R}^{1})$-valued function always has a
continuous decomposition while this does not hold for a continuous
$\mathbb{Q}_Q(\mathbb{R}^{n>1})$-valued function. See the example
in $\cite{gj}$. Nevertheless, every continuous multiple-valued
function defined on a closed interval can always be split into
continuous single-valued functions (see $\cite{gj}$ proposition
5.2).\\
Now our question is whether differentiability can be preserved
also. The example $[[x]]+[[-x]],x\in[-1,1]$ suggests that
``affinely approximatable" (see $\cite{af}$) is not enough to
guarantee a differentiable selection. Under stronger condition, we
can prove:
\begin{theorem}
If $f:[a,b]\subset\mathbb{R}\to\mathbb{Q}$ is a continuous
function and $x_0\in(a,b)$. Suppose $f$ is strongly affinely
approximatable (see $\cite{af}$) at $x_0$, then there exist
continuous functions
$f_1,f_2,\cdot\cdot\cdot,f_Q:[a,b]\to\mathbb{R}^n$ such that
$f=\sum_{i=1}^Q [[f_i]]$, and each $f_i$ is differentiable at
$x_0$.
\end{theorem}
Acknowledgement. The author would like to thank his thesis advisor
Professor Robert Hardt for introducing this subject and the paper
$\cite{gj}$. This paper would not exist without the constantly
challenging discussion with him and all the support and
encouragement he has being given in the last three years.
\section{Preliminaries}
Most of the notations, definitions and known results about
multiple-valued functions that we need can be found in
$\cite{zw1}$. The reader is also referred to $\cite{af}$ for more
details. We use standard terminology in geometric measure theory,
all of which can be found on page 669-671 of the treatise {\it Geometric Measure Theory} by H. Federer $\cite{fh}$.\\
For reader's convenience, here we state some useful results not
included in $\cite{zw1}$. The proofs of them can be found in
$\cite{af}$.
\begin{theorem}[$\cite{af},\S 1.2$]
Suppose $-\infty<r(1)\le r(2)\le\cdot\cdot\cdot\le r(Q)<\infty$
and $-\infty<s(1)\le s(2)\le \cdot\cdot\cdot \le s(Q)<\infty$,
then
$$\sum_{i=1}^Q [r(i)-s(i)]^2=\inf_{\sigma}\{\sum_{i=1}^Q [r(i)-s(\sigma(i))]^2:\sigma\;\mbox{is a permutation
of}\;\{1,2,\cdot\cdot\cdot,Q\}\}.$$
\end{theorem}
\begin{theorem}[$\cite{af},\S 2.14$]
Suppose $f:[-1,1]\rightarrow \mathbb{Q}$ such that
$f|(-1,1)\in\mathcal{Y}_2((-1,1),\mathbb{Q})$ is strictly defined
and is Dir minimizing. Then there exists
$J\in\{1,2,\cdot\cdot\cdot,Q\}$,
$k_1,k_2,\cdot\cdot\cdot,k_J\in\{1,2,\cdot\cdot\cdot,Q\}$ with
$Q=k_1+k_2+\cdot\cdot\cdot+k_J$, and
$f_1,f_2,\cdot\cdot\cdot,f_J\in \mathbb{A}(1,n)$ such that\\
(1) Whenever $-1<x<1$, and $i,j\in\{1,2,\cdot\cdot\cdot,J\}$ with
$i\not= j$, $f_i(x)\not= f_j(x)$.\\
(2) For each $-1<x<1$, $f(x)=\sum_{i=1}^J k_i[[f_i(x)]].$
\end{theorem}
\begin{theorem}[$\cite{af},\S 2.14$]
Suppose
$f_1,f_2,\cdot\cdot\cdot,f_Q\in\mathcal{Y}_2(\mathbb{B}_1^m(0),\mathbb{R}^n)$
are strictly defined, then $f=\sum_{i=1}^Q
[[f_i]]\in\mathcal{Y}_2(\mathbb{B}_1^m(0),\mathbb{Q})$.
Furthermore, in case $f$ is Dirichlet minimizing, so is each
$f_i,i=1,2,\cdot\cdot\cdot,Q$.
\end{theorem}
\section{Metric Analysis on $\mathbb{Q}$}
\subsection{$(\mathbb{Q},\mathcal{G})$ is a PC-space}
\begin{theorem} $(\mathbb{Q}, \mathcal{G})$ is a geodesic length
space, namely, given any two points $A,B\in \mathbb{Q}$, there
exists a curve $\gamma$ connecting $A$ and $B$ such that the
length of $\gamma$ equals $\mathcal{G}(A,B)$. This curve is called
a geodesic.
\end{theorem}
\begin{proof}
Given $A=\sum_{i=1}^Q [[a_i]],B=\sum_{i=1}^Q [[b_i]]$, from
Theorem 2.2 in $\cite{af}$ there is a curve
$\gamma:[0,1]\to\mathbb{Q}$ such that $\gamma(0)=A,\gamma(1)=B$,
$\gamma\in \mathcal{Y}_2((-1,1),\mathbb{Q}),$ $\gamma$ is strictly
defined and Dirichlet minimizing. We will show that this curve is
a geodesic between $A$ and $B$.\\
From Theorem 2.2, there exist $J\in\{1,2,\cdot\cdot\cdot,Q\}$,
$k_1,k_2,\cdot\cdot\cdot,k_J\in\{1,2,\cdot\cdot\cdot,Q\}$ with
$Q=k_1+k_2+\cdot\cdot\cdot+k_J$, and
$f_1,f_2,\cdot\cdot\cdot,f_J\in \mathbb{A}(1,n)$ such that
$$\gamma(t)=\sum_{i=1}^J k_i[[f_i(t)]],\;t\in[0,1],$$
and whenever $0<t<1$, and $i,j\in\{1,2,\cdot\cdot\cdot,J\}$ with
$i\not= j$, $f_i(t)\not= f_j(t)$.\\
Hence we can rewrite $\gamma$ as:
$$\gamma(t)=\sum_{i=1}^J k_i[[(1-t)p_i+tq_i]].$$
Obviously, we have
$$A=\sum_{i=1}^Q [[a_i]]=\sum_{i=1}^J
k_i[[p_i]], B=\sum_{i=1}^Q [[b_i]]=\sum_{i=1}^J k_i[[q_i]].$$
Claim 1: $\mathcal{G}^2(A,B)=\inf_\sigma\{\sum_{i=1}^Q
|a_i-b_{\sigma(i)}|^2\}=\sum_{i=1}^J k_i|p_i-q_i|^2.$\\
This is because we choose any permutation $\sigma$ of
$\{1,2,\cdot\cdot\cdot,Q\}$, and define a multiple-valued function
$f_\sigma:[0,1]\rightarrow \mathbb{Q}$ as:
$$f_\sigma(t)=\sum_{i=1}^Q [[(1-t)a_i+tb_{\sigma(i)}]].$$
Apparently, $f_\sigma(0)=A,f_\sigma(1)=B$, and
$f_\sigma\in\mathcal{Y}_2((0,1),\mathbb{Q})$, strictly defined.
Moreover
$$Dir(f_\sigma;[0,1])=\int_0^1 |Df_\sigma|^2 dt=\int_0^1
\sum_{i=1}^Q|a_i-b_{\sigma(i)}|^2 dt=\sum_{i=1}^Q
|a_i-b_{\sigma(i)}|^2.$$ Similarly,
$$Dir(\gamma;[0,1])=\sum_{i=1}^J k_i|p_i-q_i|^2.$$
Those two equalities combined with the fact that $\gamma$ is
Dirichlet minimizing proves the first claim.\\
Claim 2: $\gamma$ has constant speed, namely,
$$\mathcal{G}(\gamma(t),\gamma(s))=(t-s)\mathcal{G}(A,B),\forall
0\le s\le t\le 1.$$ Proof of Claim 2: Observe that
$\gamma|_{[s,t]}$ is also Dirichlet minimizing. Claim 1 says that
the distance between two points is realized by matching components
according to a Dirichlet minimizer connecting them. Therefore
\begin{equation*}
\begin{split}
\mathcal{G}^2(\gamma(t),\gamma(s))&=\mathcal{G}^2(\sum_{i=1}^J
k_i[[(1-t)p_i+tq_i]],\sum_{i=1}^J
k_i[[(1-s)p_i+sq_i]])\\
&=\sum_{i=1}^J k_i|(1-t)p_i+tq_i-(1-s)p_i-sq_i|^2\\
&=(t-s)^2\sum_{i=1}^J k_i|p_i-q_i|^2\\
&=(t-s)^2 \mathcal{G}^2(A,B).
\end{split}
\end{equation*}
We now can easily see that the length of $\gamma$ equals
$\mathcal{G}(A,B)$, i.e, $\gamma$ is a geodesic connecting $A$
with $B$.
\end{proof}
\begin{corollary} Using the same notations, if $n=1$, i.e, in
$\mathbb{Q}_Q(\mathbb{R})$, both of the sequences
$\{p_i\},\{q_i\}$ are non-decreasing or non-increasing. Hence
geodesic in $\mathbb{Q}_Q(\mathbb{R})$ is unique. More precisely,
two points are connected by matching components in the order of
height.
\end{corollary}
\begin{proof}
It follows directly from Theorem 2.1.
\end{proof}
\begin{remark}(1) From the proof, we also can conclude that any
geodesic connecting $A$ with $B$ must be a Dirichlet minimizer
with boundary values $A$ and $B$.\\
(2) Generally, this geodesic is not unique. For example, let
$A=[[(0,1)]]+[[(0,-1)]]\in\mathbb{Q}_2(\mathbb{R}^2),
B=[[(-1,0)]]+[[(1,0)]]\in\mathbb{Q}_2(\mathbb{R}^2).$ It is easy
to check that both
$$\gamma_1(t)=[[(-t,t-1)]]+[[(t,1-t)]],\gamma_2(t)=[[(-t,1-t)]]+[[(t,t-1)]]$$
are geodesics connecting $A$ with $B$.\\
(3) As we said, restriction of a geodesic always gives a geodesic
by Dirichlet minimality. But an extension of a geodesic may not be
a geodesic anymore. For example, considering $\gamma_1$ in the
previous one. If we extend the domain of $\gamma_1$ to $[0,2]$, it
is no longer a geodesic. In fact, the curve
$\tilde{\gamma}(t)=[[(-2t,1)]]+[[(2t,-1)]]$ is the geodesic
connecting $A$ with $[[(-2,1)]]+[[(2,-1)]]$.
\end{remark}
\begin{theorem}
The space $\mathbb{Q}$ with metric $\mathcal{G}$ is
positively-curved (a PC-space) in the sense of Alexandrov.
\end{theorem}
\begin{proof}
Given any two points $A=\sum_{i=1}^Q [[a_i]],B=\sum_{i=1}^Q
[[b_i]]\in\mathbb{Q}$, suppose $\gamma:[0,1]\rightarrow \mathbb{Q}$
is a geodesic connecting them. Take any point $C=\sum_{i=1}^Q
[[c_i]]\in\mathbb{Q}$, we will show the following inequality:
$$\mathcal{G}^2(\gamma(t),C)\ge
(1-t)\mathcal{G}^2(A,C)+t\mathcal{G}^2(B,C)-t(1-t)\mathcal{G}^2(A,B),\;\mbox{for
any}\;t\in[0,1].$$ Let's write
$$\gamma(t)=\sum_{i=1}^J k_i[[(1-t)p_i+tq_i]]=\sum_{j=1}^Q [[(1-t)p'_j+tq'_j]],$$
where in the second equality, we use the following convention:
$$p'_j=p_i,q'_j=q_i\;\mbox{if}\;k_1+\cdot\cdot\cdot+k_{i-1}+1\le j\le k_1+\cdot\cdot\cdot+k_i,$$
and $k_0=0$.\\
Fix any permutation $\sigma$ of $\{1,2,\cdot\cdot\cdot,Q\}$,
$$\sum_{j=1}^Q |(1-t)p'_j+tq'_j-c_{\sigma(j)}|^2=\sum_{j=1}^Q
|(1-t)p'_j-(1-t)c_{\sigma(j)}+tq'_j-tc_{\sigma(j)}|^2$$
\begin{equation*}
\begin{split}
&=\sum_{j=1}^Q
(1-t)^2|p'_j-c_{\sigma(j)}|^2+t^2|q'_j-c_{\sigma(j)}|^2+2t(1-t)<p'_j-c_{\sigma(j)},q'_j-c_{\sigma(j)}>\\
&=\sum_{j=1}^Q
(1-t)|p'_j-c_{\sigma(j)}|^2+t|q'_j-c_{\sigma(j)}|^2-t(1-t)|p'_j-q'_j|^2\\
&=(1-t)\sum_{j=1}^Q |p'_j-c_{\sigma(j)}|^2+t\sum_{j=1}^Q
|q'_j-c_{\sigma(j)}|^2-t(1-t)\sum_{j=1}^Q |p'_j-q'_j|^2\\
&\ge (1-t)\mathcal{G}^2(A,C)+t\mathcal{G}^2(B,C)-t(1-t)\sum_{j=1}^Q
|p'_j-q'_j|^2
\end{split}
\end{equation*}
By the definition of $p'_j,q'_j$ and the claim 1 in Theorem 3.1,
$$\sum_{j=1}^Q
|p'_j-q'_j|^2=\mathcal{G}^2(A,B).$$
Hence
$$\sum_{j=1}^Q |(1-t)p'_j+tq'_j-c_{\sigma(j)}|^2\ge
(1-t)\mathcal{G}^2(A,C)+t\mathcal{G}^2(B,C)-t(1-t)\mathcal{G}^2(A,B).$$
Taking the infimum of $\sigma$ at the left side of the above
inequality gives
$$\mathcal{G}^2(\gamma(t),C)\ge (1-t)\mathcal{G}^2(A,C)+t\mathcal{G}^2(B,C)-t(1-t)\mathcal{G}^2(A,B).$$
\end{proof}
In the special case $n=1$, we will show that in fact
$\mathbb{Q}_Q(\mathbb{R})$ is flat. Before we prove that, we need
some lemma:
\begin{lemma} Let $A=\sum_{i=1}^Q[[a_i]],B=\sum_{i=1}^Q[[b_i]]$ be any two
points in $\mathbb{Q}_Q(\mathbb{R})$. Suppose $\gamma$ is the
geodesic connecting $A$ and $B$ and can be written as:
$$\gamma(t)=\sum_{i=1}^Q [[(1-t)a_i+tb_i]].$$
Take any point $C=\sum_{i=1}^Q[[c_i]]\in \mathbb{Q}_Q(\mathbb{R})$
and suppose $\sigma$ is a permutation of
$\{1,2,\cdot\cdot\cdot,Q\}$ such that
$$\sum_{i=1}^Q(a_i-c_{\sigma(i)})^2=\mathcal{G}^2(A,C),$$ then
$$\sum_{i=1}^Q(b_i-c_{\sigma(i)})^2=\mathcal{G}^2(B,C).$$
\end{lemma}
\begin{proof}
By Corollary 3.1 we may assume that $\{a_i\}$ and $\{b_i\}$ are
both non-decreasing sequences. Applying the Theorem 2.1 to
$\mathcal{G}^2(A,C)$ we know $\{c_{\sigma(i)}\}$ is also a
non-decreasing sequence. Therefore, applying the Theorem 2.1 again
to $\mathcal{G}^2(B,C)$ gives us the desired result.
\end{proof}
\begin{theorem}
$(\mathbb{Q}_Q(\mathbb{R}),\mathcal{G})$ is a flat metric space in
the sense of Alexandrov.
\end{theorem}
\begin{proof} For any two points
$A=\sum_{i=1}^Q[[a_i]],B=\sum_{i=1}^Q[[b_i]]\in\mathbb{Q}_Q(\mathbb{R})$,
suppose
$$\gamma(t)=\sum_{i=1}^Q[[(1-t)a_i+tb_i]]$$
is the geodesic connecting $A$ with $B$. \\
From proof of Theorem 3.2, for any permutation $\sigma$ of
$\{1,2,\cdot\cdot\cdot,Q\}$,
$$\sum_{i=1}^Q |(1-t)a_i+tb_i-c_{\sigma(i)}|^2=$$
$$=(1-t)\sum_{i=1}^Q |a_i-c_{\sigma(i)}|^2+t\sum_{i=1}^Q
|b_i-c_{\sigma(i)}|^2-t(1-t)\mathcal{G}^2(A,B).$$ Now take a
permutation $\sigma$ such that $$\sum_{i=1}^Q
|a_i-c_{\sigma(i)}|^2=\mathcal{G}^2(A,C).$$ By Lemma 3.1,
$$\sum_{i=1}^Q
|b_i-c_{\sigma(i)}|^2=\mathcal{G}^2(B,C).$$ So for this
permutation $\sigma$,
$$\sum_{i=1}^Q |(1-t)a_i+tb_i-c_{\sigma(i)}|^2=(1-t)\mathcal{G}^2(A,C)+t\mathcal{G}^2(B,C)-t(1-t)\mathcal{G}^2(A,B).$$
By the definition of $\mathcal{G}^2(\gamma(t),C)$, we have
$$\mathcal{G}^2(\gamma(t),C)\le
(1-t)\mathcal{G}^2(A,C)+t\mathcal{G}^2(B,C)-t(1-t)\mathcal{G}^2(A,B).$$
The above inequality combined with Theorem 3.2 gives the following
equality:
$$\mathcal{G}^2(\gamma(t),C)=
(1-t)\mathcal{G}^2(A,C)+t\mathcal{G}^2(B,C)-t(1-t)\mathcal{G}^2(A,B),$$
which finishes the proof.
\end{proof}
\begin{remark} Generally the space $\mathbb{Q}$ is not flat, as shown by the following example.
\end{remark}
Example: Consider $\mathbb{Q}_2(\mathbb{R}^2)$,
$A=[[(0,1)]]+[[(0,0)]],B=[[(0,0)]]+[[(1,-1/2)]],$\\
$C= [[(0,0)]]+[[(-1,-1)]]$. It is easy to compute that
$$\mathcal{G}^2(A,B)=9/4,\mathcal{G}^2(A,C)=3,\mathcal{G}^2(B,C)=13/4,$$
$$\gamma(t)=[[(0,1-t)]]+[[(t,-t/2)]].$$
As for $\mathcal{G}^2(\gamma(t),C)$, there are two permutations
involved:
$$|(0,0)-(0,1-t)|^2+|(-1,-1)-(t,-t/2)|^2=9t^2/4-t+3,$$
$$|(0,0)-(t,-t/2)|^2+|(-1,-1)-(0,1-t)|^2=9t^2/4-4t+5.$$
Hence when $t\in[0,2/3]$, $\mathcal{G}^2(\gamma(t),C)=9t^2/4-t+3$
and when $t\in[2/3,1]$,
$\mathcal{G}^2(\gamma(t),C)=9t^2/4-4t+5$.\\
Now for any $t\in(0,2/3)$, \begin{equation*}
\begin{split}
\mathcal{G}^2(\gamma(t),C)&=9t^2/4-t+3\\
&>(1-t)\mathcal{G}^2(A,C)+t\mathcal{G}^2(B,C)-t(1-t)\mathcal{G}^2(A,B)
=9t^2/4-2t+3.
\end{split}
\end{equation*}
We conclude this section by a description of
$(\mathbb{Q},\mathcal{G})$ in terms of metric:
\begin{theorem}
$(\mathbb{Q},\mathcal{G})$ is a complete, separable, path
connected and locally compact metric space.
\end{theorem}
\begin{proof}
Taking a Cauchy sequence $\{A_i\}\subset\mathbb{Q}$, we consider
the sequence $\{\xi(A_i)\}\subset
\mathbb{Q}^*\subset\mathbb{R}^{PQ}$. Because
$\mbox{Lip}(\xi)<\infty$, we infer
$$\{\xi(A_i)\}\;\mbox{is a Cauchy sequence in}\;\mathbb{R}^{PQ}.$$
Due to the completeness of $\mathbb{R}^{PQ}$, there is an element
$A\in\mathbb{R}^{PQ}$, such that
$$\lim_{i\to\infty}\xi(A_i)=A.$$
Because $\mathbb{Q}^*$ is closed in $\mathbb{R}^{PQ}$, we conclude
that $A\in\mathbb{Q}^*$. Therefore $\lim_{i\to\infty}
A_i=\xi^{-1}(A)$ thanks to the fact that
$\mbox{Lip}(\xi^{-1})<\infty$.\\
As for the separability, one can check that the set
$$\{\sum_{i=1}^Q[[a_i]],a_i\;\mbox{is a rational point
in}\;\mathbb{R}^n\}$$ is a countable dense subset of
$\mathbb{Q}$.\\
Take any two points $A,B\in\mathbb{Q}$, any geodesic connecting
them gives a path between them. Hence $\mathbb{Q}$ is path
connected.\\
As for the locally compactness, it is obviously once we notice the
bi-Lipschitzian correspondence between $\mathbb{Q}$ and
$\mathbb{Q}^*\subset\mathbb{R}^{PQ}$.
\end{proof}
\begin{remark}
By Hopf-Rinow Theorem, any bounded closed set in $\mathbb{Q}$ is
compact.
\end{remark}
\subsection{Metric Tangent Cone of $(\mathbb{Q},\mathcal{G})$}
\begin{definition}
For $A,B,C\in\mathbb{Q}$, define
$$\alpha(A;B,C)=\DF{\mathcal{G}^2(A,B)+\mathcal{G}^2(A,C)-\mathcal{G}^2(B,C)}{2\mathcal{G}(A,B)\mathcal{G}(A,C)},A\not= B,C.$$
\end{definition}
\begin{theorem}
Let $\gamma_1,\gamma_2$ be geodesics starting from $A$, then the
function
$$t,s\in(0,1]\to\alpha(A;\gamma_1(t),\gamma_2(s))\;\mbox{is
nondecreasing in}\;s,t.$$ The angle
$\angle(\gamma_1,\gamma_2)\in[0,\pi]$ between $\gamma_1$ and
$\gamma_2$ is thus defined by the formula
$$\cos(\angle(\gamma_1,\gamma_2)):=\inf_{s,t}\alpha(A;\gamma_1(t),\gamma_2(s))=\lim_{s,t\downarrow
0}\alpha(A;\gamma_1(t),\gamma_2(s)).$$
\end{theorem}
\begin{proof}
See $\cite{ags}$ Lemma 12.3.4.
\end{proof}
For a fixed $A\in\mathbb{Q}$ let us denote by $G(A)$ the set of
all geodesics $\gamma$ starting from $A$ and parameterized in some
interval $[0,T_{\gamma}]$; recall that the metric velocity of
$\gamma$ is $|\gamma'|=\mathcal{G}(\gamma(t),A)/t,t\in(0,T]$. We
set
$$||\gamma||_A:=|\gamma'|,<\gamma_1,\gamma_2>_A:=||\gamma_1||_A||\gamma_2||_A\cos(\angle(\gamma_1,\gamma_2)),$$
$$d_A^2(\gamma_1,\gamma_2):=||\gamma_1||_A^2+||\gamma_2||_A^2-2<\gamma_1,\gamma_2>_A.$$
If $\gamma\in G(A)$ and $\lambda>0$ we denote by $\lambda\gamma$
the geodesic
$$(\lambda\gamma)_t:=\gamma_{\lambda
t},T_{\lambda\gamma}=\lambda^{-1}T_{\gamma},$$ and we observe that
for each $\gamma_1,\gamma_2\in G(A),\lambda>0,$ it holds
$$||\lambda\gamma||_A=\lambda||\gamma||_A,<\lambda\gamma_1,\gamma_2>_A=<\gamma_1,\lambda\gamma_2>_A=\lambda<\gamma_1,\gamma_2>_A.$$
Recall that the restriction of a geodesic is still a geodesic; we
say that $\gamma_1\sim\gamma_2$ if there exists $\epsilon>0$ such
that $\gamma_1|_{[0,\epsilon]}=\gamma_2|_{[0,\epsilon]}$.
\begin{theorem}
If $\gamma_1,\gamma_2:[0,T]\to\mathbb{Q}$ are two geodesics
starting from $A$, we have
$$d_A(\gamma_1,\gamma_2)=\lim_{t\downarrow
0}\DF{\mathcal{G}(\gamma_1(t),\gamma_2(t))}{t}.$$ In particular,
the function $d_A$ defined above is a distance on the quotient
space $G(A)/\sim$. The completion of $G(A)/\sim$ is called the
tangent cone $\bold{Tan}_A(\mathbb{Q})$ at the point $A$.
\end{theorem}
\begin{proof}
It follows from the same argument as $\cite{ags}$ Theorem 12.3.6.
\end{proof}
Before we give an explicit representation of the abstract tangent
cone $\bold{Tan}_A(\mathbb{Q})$, we give several definitions
concerning representation of elements in $\mathbb{Q}$.
\begin{definition} Define
$\sigma:\mathbb{Q}\to\{1,2,\cdot\cdot\cdot,Q\}$ by
$$\sigma(x)=\mbox{card}[\mbox{spt}(x)].$$
\end{definition}
\begin{remark}
It is easy to see that $\sigma$ is lower semi-continuous.
\end{remark}
\begin{definition}
For any $x\in\mathbb{Q}$, which can be written as
$$x=\sum_{i=1}^J k_i[[x_i]],$$
for some
$J\in\{1,2,\cdot\cdot\cdot,Q\},k_i\in\{1,2,\cdot\cdot\cdot,Q\},k_1\le
k_2\le\cdot\cdot\cdot\le k_J, \sum_{i=1}^J k_i=Q$ and $x_i$
distinct points in $\mathbb{R}^n$, define the signature of $x$ as:
$$S(x)=(J,k_1,k_2,\cdot\cdot\cdot,k_J).$$
\end{definition}
\begin{remark}
In $S(x)$, $J$ is exactly $\sigma(x)$.
\end{remark}
\begin{definition}
For a fixed positive integer $Q$,
$(J,k_1,k_2,\cdot\cdot\cdot,k_J)$ is called a permissible
decomposition of $Q$ if
$$J\in\{1,2,\cdot\cdot\cdot,Q\},k_i\in\{1,2,\cdot\cdot\cdot,Q\},k_1\le
k_2\le\cdot\cdot\cdot\le k_J,\sum_{i=1}^J k_i=Q.$$ The set of all
permissible decompositions of $Q$ is denoted as $\mathcal{P}(Q)$.
\end{definition}
\begin{remark}
From $\cite{lw}$ (15.8), $\mbox{card}(\mathcal{P}(Q))=p(Q)$, where
$p(n)$ is defined to be the number of unordered partitions of $n$.
\end{remark}
\begin{proposition}
Suppose $S(A)=(J,k_1,\cdot\cdot\cdot,k_J)$. Let
$\gamma_1,\gamma_2\in\ G(A)$ such that
\begin{equation*}
\begin{split}
\gamma_1(t)&=\sum_{i=1}^J\sum_{j=k_1+\cdot\cdot\cdot+k_{i-1}+1}^{k_1+\cdot\cdot\cdot+k_i}
[[a_i+tv_j]],\\
\gamma_2(t)&=\sum_{i=1}^J\sum_{j=k_1+\cdot\cdot\cdot+k_{i-1}+1}^{k_1+\cdot\cdot\cdot+k_i}
[[a_i+tw_j]].
\end{split}
\end{equation*}
Then
$$d_A^2(\gamma_1,\gamma_2)=\sum_{i=1}^J
\mathcal{G}^2(\sum_{j=k_1+\cdot\cdot\cdot+k_{i-1}+1}^{k_1+
\cdot\cdot\cdot+k_i}[[v_j]],\sum_{j=k_1+\cdot\cdot\cdot+k_{i-1}+1}^{k_1+
\cdot\cdot\cdot+k_i}[[w_j]]).$$
\end{proposition}
\begin{proof}
When $t$ small enough, the distance between $\gamma_1(t)$ and
$\gamma_2(t)$ is obtained by: for each
$i\in\{1,2,\cdot\cdot\cdot,J\}$,
$$\sum_{j=k_1+\cdot\cdot\cdot+k_{i-1}+1}^{k_1+
\cdot\cdot\cdot+k_i}[[a_i+tv_j]]\;\mbox{matching
with}\;\sum_{j=k_1+\cdot\cdot\cdot+k_{i-1}+1}^{k_1+
\cdot\cdot\cdot+k_i}[[a_i+tw_j]].$$ Hence using Theorem 3.6
finishes the proof.
\end{proof}
This proposition gives us the description of
$\bold{Tan}_A(\mathbb{Q})$:
\begin{theorem} For any point $A\in\mathbb{Q}$, such that
$S(A)=(J,k_1,\cdot\cdot\cdot,k_J)$, we have the isometry
$$\bold{Tan}_A(\mathbb{Q})\cong\mathbb{Q}_{k_1}(\mathbb{R}^n)\times
\mathbb{Q}_{k_2}(\mathbb{R}^n)\times \cdot\cdot\cdot\times
\mathbb{Q}_{k_J}(\mathbb{R}^n)$$ with the product metric.
\end{theorem}
\begin{proof}
Observe that any geodesic starting with $A$ is uniquely determined
by the initial velocity, which is given by
$$(\sum_{i=1}^{k_1}[[v_i]],\sum_{i=k_1+1}^{k_1+k_2}[[v_i]],\cdot\cdot\cdot,
\sum_{k_1+\cdot\cdot\cdot+k_{J-1}+1}^Q[[v_i]]).$$ Proposition 3.1
shows the metric on $\bold{Tan}_A(\mathbb{Q})$. We are done.
\end{proof}
\begin{remark}
For any two points $x,y\in\mathbb{Q}$,
$$\bold{Tan}_x(\mathbb{Q})\cong\bold{Tan}_y(\mathbb{Q})\iff S(x)=S(y).$$
\end{remark}
\subsection{Exponential Map}
Given a point $A=\sum_{i=1}^J k_i[[a_i]]$ of signature
$S(A)=(J,k_1,\cdot\cdot\cdot,k_J)$ and a nonzero element
$\vec{v}=(\sum_{i=1}^{k_1}
[[v_i]],\cdot\cdot\cdot,\sum_{i=k_1+\cdot\cdot\cdot+k_{J-1}+1}^Q[[v_i]])\in
\bold{Tan}_A(\mathbb{Q})$, there is a unique parameterized
geodesic $\gamma:(-\epsilon,\epsilon)\to\mathbb{Q}$ of the form
$$\gamma(t)=\sum_{i=1}^J\sum_{j=k_1+\cdot\cdot\cdot+k_{i-1}+1}^{k_1+\cdot\cdot\cdot+k_i}[[a_i+tv_j]].$$
To indicate the dependence of this geodesic on the vector
$\vec{v}$, we denote it by $\gamma(t,\vec{v})=\gamma$. Obviously
we have
\begin{proposition}
If the geodesic $\gamma(t,\vec{v})$ is defined for
$t\in(-\epsilon,\epsilon)$, then the geodesic
$\gamma(t,\lambda\vec{v}),\lambda\in\mathbb{R},\lambda>0,$ is
defined for $t\in(-\epsilon/\lambda,\epsilon/\lambda),$ and
$\gamma(t,\lambda\vec{v})=\gamma(\lambda t,\vec{v})$.
\end{proposition}
\begin{definition}
If $\vec{v}\in\bold{Tan}_A(\mathbb{Q}),\vec{v}\not=\vec{0}$, is
such that $\gamma(|\vec{v}|,\vec{v}/|\vec{v}|)=\gamma(1,\vec{v})$
is defined, we set
$$\exp_A(\vec{v})=\gamma(1,\vec{v})\;\mbox{and}\;\exp_A(\vec{0})=A.$$
\end{definition}
\begin{remark}
Exponential map is not necessarily one-to-one, as shown in Remark
3.1(1). However, if restricted in a small neighborhood of $0$ in
$\bold{Tan}_A(\mathbb{Q})$, it is an isometry.
\end{remark}
\begin{theorem}
For any $A\in\mathbb{Q}$ with signature
$S(A)=(J,k_1,\cdot\cdot\cdot,k_J)$, there is a positive $\epsilon$
such that $\exp_A:\mathbb{B}_\epsilon(0)\subset
\bold{Tan}_A(\mathbb{Q})\to\mathbb{Q}$ is an isometry.
\end{theorem}
\begin{proof}
Suppose $A=\sum_{i=1}^J k_i[[a_i]]$, and let
$$\delta=2^{-1}\inf\{|a_i-a_j|,i\not=j\}.$$
Choose $\epsilon=\delta/2$. Let $$\vec{v}=(\sum_{i=1}^{k_1}
[[v_i]],\cdot\cdot\cdot,\sum_{i=k_1+\cdot\cdot\cdot+k_{J-1}+1}^Q[[v_i]])\in
\mathbb{B}_\epsilon(0)\subset \bold{Tan}_A(\mathbb{Q}),$$
$$\vec{w}=(\sum_{i=1}^{k_1}
[[w_i]],\cdot\cdot\cdot,\sum_{i=k_1+\cdot\cdot\cdot+k_{J-1}+1}^Q[[w_i]])\in\mathbb{B}_\epsilon(0)\subset
\bold{Tan}_A(\mathbb{Q}).$$ Therefore,
$$\gamma_1(t)=\sum_{i=1}^J\sum_{j=k_1+\cdot\cdot\cdot+k_{i-1}+1}^{k_1+\cdot\cdot\cdot+k_i}[[a_i+tv_j]],$$
$$\gamma_2(t)=\sum_{i=1}^J\sum_{j=k_1+\cdot\cdot\cdot+k_{i-1}+1}^{k_1+\cdot\cdot\cdot+k_i}[[a_i+tw_j]],$$
$$\exp_A(\vec{v})=\sum_{i=1}^J\sum_{j=k_1+\cdot\cdot\cdot+k_{i-1}+1}^{k_1+\cdot\cdot\cdot+k_i}[[a_i+v_j]],$$
$$\exp_A(\vec{w})=\sum_{i=1}^J\sum_{j=k_1+\cdot\cdot\cdot+k_{i-1}+1}^{k_1+\cdot\cdot\cdot+k_i}[[a_i+w_j]].$$
Then
\begin{equation*}
\begin{split}
\mathcal{G}^2(\exp_A(\vec{v}),\exp_A(\vec{w}))&=d^2(\vec{v},\vec{w})\\
&=\sum_{i=1}^J
\mathcal{G}^2(\sum_{j=k_1+\cdot\cdot\cdot+k_{i-1}+1}^{k_1+
\cdot\cdot\cdot+k_i}[[v_j]],\sum_{j=k_1+\cdot\cdot\cdot+k_{i-1}+1}^{k_1+
\cdot\cdot\cdot+k_i}[[w_j]]).
\end{split}
\end{equation*}
\end{proof}
\subsection{A Decomposition of $(\mathbb{Q},\mathcal{G})$}
We decompose $\mathbb{Q}$ according to the signature:
\begin{definition}
Define
$$I_{J,k_1,k_2,\cdot\cdot\cdot,k_J}=\{x\in\mathbb{Q}:S(x)=(J,k_1,k_2,\cdot\cdot\cdot,k_J)\},$$
for any $(J,k_1,k_2,\cdot\cdot\cdot,k_J)\in\mathcal{P}(Q)$.
\end{definition}
\begin{definition}
Define
$$\mathcal{I}_i=\{x\in\mathbb{Q}:\mbox{card}(\mbox{spt}(x))=i\}=\cup_{(i,k_1,\cdot\cdot\cdot,k_i)\in\mathcal{P}(Q)} I_{i,k_1,\cdot\cdot\cdot,k_i},$$
for $i\in\{1,2,\cdot\cdot\cdot,Q\}.$
\end{definition}
\begin{remark}
Geometrically, a point in $\mathcal{I}_1$ is like the vertex point
of a cone in the sense that the tangent cone at points in
$\mathcal{I}_1$ is isometric to $\mathbb{Q}$ itself. Points in
$\mathcal{I}_Q$ are like faces of a polyhedra cone in the sense
that the tangent cone at points in $\mathcal{I}_Q$ is isometric to
$\mathbb{R}^{nQ}$, namely, locally $\mathcal{I}_Q$ is flat. We
will rigorously prove the flatness of $\mathcal{I}_Q$ in the sense
of Alexandrov later.
\end{remark}
\begin{theorem}
(1) $\mathcal{I}_i$ only can be approximated by elements in
$\mathcal{I}_j$ for $j\ge i$.\\
(2) $\mathcal{I}_Q$ is open in $\mathbb{Q}$.
\end{theorem}
\begin{proof}
They both follow from the lower semi-continuity of the function
$\sigma$.
\end{proof}
\begin{theorem}
$\mathcal{I}_Q$ is path connected and dense in $\mathbb{Q}$.
\end{theorem}
\begin{proof}
We endow $\mathbb{R}^n$ with the lexicographical order. Take any
two points $A,B\in\mathcal{I}_Q$. Write them as
$$A=\sum_{i=1}^Q[[a_i]],B=\sum_{i=1}^Q [[b_i]],$$
with $a_1<a_2<\cdot\cdot\cdot<a_Q$, and
$b_1<b_2<\cdot\cdot\cdot<b_Q$. Define a curve
$\gamma:[0,1]\to\mathbb{Q}$ as
$$\gamma(t)=\sum_{i=1}^Q [[(1-t)a_i+tb_i]].$$
Obviously $\gamma$ connects $A$ with $B$. We will show that
$\gamma\in\mathcal{I}_Q$.\\
Suppose not, i.e, there are $i<j$ and $t\in(0,1)$ such that
$$(1-t)a_i+tb_i=(1-t)a_j+tb_j.$$
Hence $a_i-a_j=\DF{t}{1-t}(b_j-b_i)$, which is impossible because
$a_i<a_j,b_i<b_j$.\\
$\mathcal{I}_Q$ is dense in $\mathbb{Q}$ because the set
$$\{\sum_{i=1}^Q[[a_i]]:a_i\;\mbox{is a rational point in}\;\mathbb{R}^n,\;\mbox{and they are distinct}\}$$
is dense in $\mathbb{Q}$.
\end{proof}
\begin{theorem} $\mathcal{I}_1$ is path connected.
\end{theorem}
\begin{proof}
Take any two points $A,B\in\mathcal{I}_1$. Write them as
$$A=Q[[a]],B=Q[[b]].$$
Then the curve
$\gamma:[0,1]\to\mathcal{I}_1,\gamma(t):=Q[[(1-t)a+tb]]$ certainly
works.
\end{proof}
\begin{remark}
Generally $\mathcal{I}_i$ is not path connected for $i\not=1,Q$.
For example, we consider $\mathbb{Q}_3(\mathbb{R})$. Let
$A=[[1]]+2[[2]],B=[[2]]+2[[1]]\in\mathcal{I}_2$. Suppose there is
a continuous curve $\gamma:[0,1]\to\mathcal{I}_2$ connecting $A$
with $B$. Since $\gamma\in\mathcal{I}_2$, we can write $\gamma$ as
$$\gamma=[[\gamma_1]]+2[[\gamma_2]],$$
for $\gamma_i:[0,1]\to\mathbb{R},i=1,2$ and
$\gamma_1(t)\not=\gamma_2(t),\forall t\in[0,1]$. Because
$\gamma_1$ and $\gamma_2$ never cross, they both are continuous.
Considering the initial conditions, we have
$$\gamma_1(0)=1,\gamma_1(1)=2,$$
$$\gamma_2(0)=2,\gamma_2(1)=1.$$
Therefore $\gamma_1$ and $\gamma_2$ must meet at some point in
$(0,1)$ due to elementary facts about continuous functions. This
shows that $\mathcal{I}_2$ is not path connected.
\end{remark}
However, we are able to prove the following theorem:
\begin{theorem} Take
$(J,k_1,k_2,\cdot\cdot\cdot,k_J)\in\mathcal{P}(Q)$ and $A\in
I_{J,k_1,k_2,\cdot\cdot\cdot,k_J}$. If $B\in
I_{J,k_1,k_2,\cdot\cdot\cdot,k_J}$, and $\mathcal{G}(A,B)$ is
small enough, then any geodesic $\gamma$ connecting them lies in
$I_{J,k_1,k_2,\cdot\cdot\cdot,k_J}$.
\end{theorem}
\begin{proof}
Let
$$A=\sum_{i=1}^J k_i[[a_i]],B=\sum_{i=1}^J k_i[[b_i]],$$
for distinct $a_i's$ and $b_i's$. Let
$$\delta=2^{-1}\inf\{|a_i-a_j|,i\not=j\}.$$
Choose $B$ close enough with $A$ such that
$$|b_i-a_i|<\delta/{4^Q},i\in\{1,2,\cdot\cdot\cdot,Q\}.$$
Obviously, $\mathcal{G}^2(A,B)=\sum_{i=1}^J k_i|a_i-b_i|^2$.
Therefore, according to proof of Theorem 3.1, any geodesic
$\gamma:[0,1]\to\mathbb{Q}$ connecting $A$ with $B$ can be written
as:
$$\gamma(t)=\sum_{i=1}^J k_i[[(1-t)a_i+tb_i]].$$
We will show that $\gamma\in I_{J,k_1,k_2,\cdot\cdot\cdot,k_J}$.
If not, i.e, there are $i<j$ and $t\in(0,1)$ such that
$$(1-t)a_i+tb_i=(1-t)a_j+tb_j.$$
Subtract $a_i$ from both sides,
$$t(b_i-a_i)=t(b_j-a_j)+a_j-a_i.$$
The absolute value of LHS is less than $\delta/{4^Q}$, while the
absolute value of RHS is great than
$\delta-\delta/{4^Q}=(1-4^{-Q})\delta$. A contradiction.
\end{proof}
As we promised, we will show $\mathcal{I}_Q$ is locally flat:
\begin{theorem}
$\mathcal{I}_Q$ is locally flat in the sense of Alexandrov.
\end{theorem}
\begin{proof}
Take any point $A=\sum_{i=1}^Q\in\mathcal{I}_Q$. Let
$$\delta=2^{-1}\inf\{|a_i-a_j|,i\not=j\},$$
and
$$U=\{x\in\mathbb{Q}:\mathcal{G}(x,A)<\delta/2\}\subset\mathcal{I}_Q.$$
Take any point $B\in U$ and $\gamma:[0,1]\to\mathcal{I}_Q$ be a
geodesic connecting $A$ with $B$ (Theorem 3.12 guarantees that
$\gamma\in\mathcal{I}_Q$ once $U$ is small enough). We write
$\gamma$ as
$$\gamma(t)=\sum_{i=1}^Q[[(1-t)a_i+tb_i]],$$
which means
$$\mathcal{G}^2(A,B)=\sum_{i=1}^Q|a_i-b_i|^2.$$
Take $C\in U$ and suppose
$$\mathcal{G}^2(A,C)=\sum_{i=1}^Q|a_i-c_i|^2,$$
which guarantees that
$$\mathcal{G}^2(B,C)=\sum_{i=1}^Q|b_i-c_i|^2,\mathcal{G}^2(\gamma(t),C)=\sum_{i=1}^Q|(1-t)a_i+tb_i-c_i|^2,\forall
t.$$ Then we have the following equality
$$\mathcal{G}^2(\gamma(t),C)=(1-t)\mathcal{G}^2(A,C)+t\mathcal{G}^2(B,C)-t(1-t)\mathcal{G}^2(A,B),\forall
t\in[0,1].$$
\end{proof}
The proof suggests a one-to-one correspondence between a
neighborhood of $A$ with an open set in $\mathbb{R}^{nQ}$.
Moreover, this correspondence turns out to be an isometry:
\begin{theorem}
For any point $A\in\mathcal{I}_Q$, there is an open neighborhood
$U$ of $A$ in $\mathcal{I}_Q$ and an bijective isometry between
$U$ with some open set in $\mathbb{R}^{nQ}$.
\end{theorem}
\begin{proof}
Using the same notations as Theorem 3.13, for any $B,C\in U$, we
have
$$\mathcal{G}^2(A,B)=\sum_{i=1}^Q|a_i-b_i|^2,\mathcal{G}^2(A,C)=\sum_{i=1}^Q|a_i-c_i|^2,$$
$$\mathcal{G}^2(B,C)=\sum_{i=1}^Q|b_i-c_i|^2.$$
That means that the distance between any two points in $U$ is
obtained in a unique way. Define
$V=\{x\in\mathbb{R}^{nQ}:|x-(a_1,a_2,\cdot\cdot\cdot,a_Q)|<\delta/2\}\subset\mathbb{R}^{nQ}$
and define
$$\phi:U\to
V,\phi(\sum_{i=1}^Q[[b_i]])=(b_1,b_2,\cdot\cdot\cdot,b_Q).$$ It is
easy to check that is is a bijection and an isometry.
\end{proof}
\begin{remark}
We can use the lexicographic order to relate $\mathbb{Q}$ with
$\mathbb{R}^{nQ}$. But it is not very useful because is is not
necessarily an isometry. For example, we endow $\mathbb{R}^n$ with
the lexicographical order. Then take any
$A=\sum_{i=1}^Q[[a_i]]\in\mathbb{Q}$, we can order $a_i$ as
following
$$a_1\le a_2\le \cdot\cdot\cdot\le a_Q$$
in a unique way. Hence we are able to define
$$\psi:\mathbb{Q}\to\mathbb{R}^{nQ},\psi(\sum_{i=1}^Q[[a_i]])=(a_1,a_2,\cdot\cdot\cdot,a_Q).$$
Unless $n=1$, $\psi$ will not be an isometry. For example, we
consider $\mathbb{Q}_2(\mathbb{R}^2)$. For $0<\epsilon<1/2$, let
$A_\epsilon=[[(1,1)]]+[[(1+\epsilon,2)]],B=[[(1+\epsilon,1)]]+[[(1,2)]]$.
Under the map $\psi$, we have
$$\psi(A_\epsilon)=(1,1,1+\epsilon,2),\psi(B_\epsilon)=(1,2,1+\epsilon,1).$$
Therefore,
$\mathcal{G}(A_\epsilon,B_\epsilon)=\sqrt{2}\epsilon,|\psi(A_\epsilon)-\psi(B_\epsilon)|=\sqrt{2}.$
Letting $\epsilon\downarrow 0$ shows that $\psi$ is not even a
Lipschitz map.
\end{remark}
\section{Tensor Sum of Multiple-Valued Functions}
\begin{definition}
Suppose $f(x)=\sum_{i=1}^p [[f_i(x)]]$, $g(x)=\sum_{j=1}^q
[[g_j(x)]]$, where $p$ and $q$ are not necessarily the same.
Define
$$(f\oplus g)(x)=\sum_{i,j} [[f_i(x)+g_j(x)]].$$
(i.e the Tensor sum is a $pq-$valued function). \end{definition}
\begin{remark} It is easy to check:
$$\eta(f\oplus g)=\eta(f)+\eta(g)$$
\end{remark}
\begin{remark} If $f$ and $g$ are both Dirichlet minimizing, although
their tensor sum $f\oplus g$ may not be minimizing, the average
function $\eta(f\oplus g)$ is harmonic. This is because both
$\eta(f)$ and $\eta(g)$ are harmonic. To see that $f\oplus g$ may
not be minimizing, let $f(x)=[[x]]+[[-1]],g(x)=[[1-x]]+[[-1]]$.
Their domains are both $[0,1]$. Both of the functions are
Dirichlet minimizers. But their tensor sum $f\oplus
g=[[x-1]]+[[-x]]+[[1]]+[[-2]]$ is no longer minimizing since there
is a branch point for the function $f\oplus g$.
\end{remark}
\begin{theorem}[Weighted
Triangular Inequality] Suppose $f(x)=\sum_{i=1}^p [[f_i(x)]]$,
$g(x)=\sum_{j=1}^q [[g_j(x)]]$, then
$$\mathcal{G}(f\oplus g,pq[[0]])\le q^{1/2}\mathcal{G}(f,p[[0]])+p^{1/2}\mathcal{G}(g,q[[0]]).$$
\end{theorem}
\begin{proof}
By definition,
\begin{equation*}
\begin{split}
\mathcal{G}(f\oplus g,pq[[0]])^2&=\sum_{i=1}^p \sum_{j=1}^q
|f_i(x)+g_j(x)|^2\\
&=\sum_{i=1}^p \sum_{j=1}^q (|f_i(x)|^2+
|g_j(x)|^2+2f_i(x)\cdot g_j(x))\\
&=q\mathcal{G}(f,p[[0]])^2+p\mathcal{G}(g,q[[0]])^2+2\sum_{i=1}^p
\sum_{j=1}^q f_i(x)\cdot g_j(x)\\
&=q\mathcal{G}(f,p[[0]])^2+p\mathcal{G}(g,q[[0]])^2+2pq\eta(f)\cdot
\eta(g)
\end{split}
\end{equation*}
Since for each $S,T\in \mathbb{Q}$, $|\eta(S)-\eta(T)|\le
Q^{-1/2}\mathcal{G}(S,T),$
$$|\eta(f)|\le
p^{-1/2}\mathcal{G}(f,p[[0]]),|\eta(g)|\le
q^{-1/2}\mathcal{G}(g,q[[0]]).$$ Therefore,
\begin{equation*}
\begin{split}
\mathcal{G}(f\oplus g,pq[[0]])^2&\le
q\mathcal{G}(f,p[[0]])^2+p\mathcal{G}(g,q[[0]])^2\\
&+2p\cdot q\cdot p^{-1/2} \mathcal{G}(f,p[[0]])\cdot
q^{-1/2}\mathcal{G}(g,q[[0]])\\
&=[q^{1/2}\mathcal{G}(f,p[[0]])+p^{1/2}\mathcal{G}(g,q[[0]])]^2.
\end{split}
\end{equation*}
\end{proof}
\begin{theorem} Suppose $f(x)=\sum_{i=1}^p [[f_i(x)]]$, $g(x)=\sum_{j=1}^q
[[g_j(x)]]$, where
$f_i,g_j\in\mathcal{Y}_2(\mathbb{B}_1^m(0),\mathbb{R}^n)$ are
strictly defined, then
$$\mbox{Dir}(f\oplus g)=q\mbox{Dir}(f)+p\mbox{Dir}(g)+2pq\int <ap
D(\eta\circ f),ap D(\eta\circ g)>d\mathcal{H}^m.$$
\end{theorem}
\begin{proof} By Theorem 2.3, $f\oplus
g\in\mathcal{Y}_2(\mathbb{B}_1^m(0),\mathbb{Q}_{pq}(\mathbb{R}^n))$.
It is easy to see that Dir($f\oplus g)=\sum_{j=1}^q \mbox{Dir}(f\oplus g_j)$.\\
From $\cite{af} \S 2.3$, we have
$$\mbox{Dir}(f\oplus g_j)=\mbox{Dir}(f)+2p\int <\mbox{ap D}(\eta\circ f),\mbox{ap D}(g_j)>d\mathcal{H}^m+p\mbox{Dir}(g_j).$$
Sum them up, we get
\begin{equation*}
\begin{split}
\mbox{Dir}(f\oplus g)&=q\mbox{Dir}(f)+\sum_{j=1}^q 2p\int
<\mbox{ap
D}(\eta\circ f),\mbox{ap D}(g_j)>d\mathcal{H}^m +p \mbox{Dir}(g)\\
&=q\mbox{Dir}(f)+2p\int <\mbox{ap D}(\eta\circ f),\mbox{ap
D}(\sum_{j=1}^q g_j)>d\mathcal{H}^m +p\mbox{Dir}(g)\\
&=q\mbox{Dir}(f)+p\mbox{Dir}(g)+2pq\int <\mbox{ap D}(\eta\circ
f),\mbox{ap D}(\eta\circ g)>d\mathcal{H}^m.
\end{split}
\end{equation*}
\end{proof}
\begin{definition}
For a Q-valued function $f$, its $L^k$ norm ($0<k<\infty$) is
defined to be
$$\|f\|_k=(\int \mathcal{G}(f(x),Q[[0]])^k d\mathcal{H}^m)^{1/k}.$$
Moreover, its $L^\infty$ norm is defined to be
$$\|f\|_\infty=\inf_M\{\mathcal{G}(f(x),Q[[0]])\le M, a.e\}.$$
\end{definition} \begin{theorem}[
Weighted Minkowski Inequality] Suppose $f$ is a $p$-valued
function and $g$ is a $q$-valued function, $1\le k\le \infty$,
then
$$\|f\oplus g\|_k\le q^{1/2}\|f\|_k+p^{1/2}\|g\|_k$$
\end{theorem}
\begin{proof} When $k=1,\;\mbox{or}\;\infty$, it follows easily from the Weighted Triangular Inequality.  \\
When $1<k<\infty$,
$$\int \mathcal{G}(f\oplus g,pq[[0]])^k d\mathcal{H}^m=\int \mathcal{G}(f\oplus g,pq[[0]])^{k-1}\mathcal{G}(f\oplus g,pq[[0]]) d\mathcal{H}^m$$
$$\le\int \mathcal{G}(f\oplus g,pq[[0]])^{k-1} q^{1/2}\mathcal{G}(f,p[[0]]) d\mathcal{H}^m+\int \mathcal{G}(f\oplus g,pq[[0]])^{k-1} p^{1/2}
\mathcal{G}(g,q[[0]]) d\mathcal{H}^m.$$ Now applying the
H$\ddot{o}$lder inequality for parameters $\DF{k}{k-1},k$:
$$\int \mathcal{G}(f\oplus g,pq[[0]])^k d\mathcal{H}^m\le (\int \mathcal{G}(f\oplus g,pq[[0]])^k)^{\DF{k-1}{k}}(\int (q^{1/2}\mathcal{G}(f,p[[0]]))^k)^{1/k}$$
$$+(\int \mathcal{G}(f\oplus g,pq[[0]])^k)^{\DF{k-1}{k}}(\int (p^{1/2}\mathcal{G}(g,q[[0]]))^k)^{1/k}$$
$$=\|f\oplus g\|_k^{k-1}q^{1/2}\|f\|_k+\|f\oplus g\|_k^{k-1}p^{1/2}\|g\|_k$$
If $\|f\oplus g\|_k$ is zero, we are done. Otherwise, divide both sides of the above inequality by $\|f\oplus g\|_k^{k-1}$.\\
\end{proof}
\begin{theorem} Suppose $f$ is a $p$-valued
function and $g$ is a $q$-valued function, $0<k<1$, then
$$\|f\oplus g\|_k^k\le q^{k/2}\|f\|_k^k+p^{k/2}\|g\|_k^k.$$
\end{theorem}
\begin{proof} Using the inequality $(a+b)^k\le a^k+b^k, a\ge 0,b\ge 0,
0<k<1$ and the Weighted Triangular Inequality, we have
\begin{equation*}
\begin{split}
\mathcal{G}(f\oplus g,pq[[0]])^k&\le
(q^{1/2}\mathcal{G}(f,p[[0]])+p^{1/2}\mathcal{G}(g,q[[0]]))^k\\
&\le[q^{1/2}
\mathcal{G}(f,p[[0]])]^k+[p^{1/2}\mathcal{G}(g,q[[0]])]^k.
\end{split}
\end{equation*}
Hence \begin{equation*} \begin{split} \|f\oplus g\|_k^k&=\int
\mathcal{G}(f\oplus g,pq[[0]])^k d\mathcal{H}^m\\ &\le \int
q^{k/2}\mathcal{G}(f,p[[0]])^k d\mathcal{H}^m+ \int p^{k/2}
\mathcal{G}(g,q[[0]])^k d\mathcal{H}^m\\
&=q^{k/2}\|f\|_k^k+p^{k/2}\|g\|_k^k.
\end{split}
\end{equation*}
\end{proof}
\section{Derivative of Multiple-Valued Functions}
\subsection{Definition of Derivative}
For reader's convenience, we include here the definition of
derivative used in $\cite{af}$.
\begin{definition}
(a) $f$ is called affine if there are $A_1,\cdot\cdot\cdot,A_Q$
where each $A_i$ is an affine map from $\mathbb{R}^m$ to
$\mathbb{R}^n$, such that
$$f(x)=\sum_{i=1}^Q [[A_i(x)]].$$
(b) $f$ is called affinely approximatable at $x_0$ if there are
affine maps $A_1,\cdot\cdot\cdot,A_Q$ from $\mathbb{R}^m$ to
$\mathbb{R}^n$ such that
$$\lim_{|x-x_0|\to 0}\DF{\mathcal{G}(f(x),\sum_{i=1}^Q [[A_i(x)]])}{|x-x_0|}=0.$$
(c) $f$ is strongly affinely approximatable at $x_0$ if (b) holds
for $f$ at $x_0$ and $A_i=A_j$ if $A_i(x_0)=A_j(x_0)$.
\end{definition}
\begin{remark}
If $f$ is affinely approximatable at $x_0$ with
$\sum_{i=1}^Q[[A_i]]$ as its affine approximation, then obviously
$f(x_0)= \sum_{i=1}^Q [[A_i(x_0)]]$ and
$A_i(x)=A_i(x_0)+L_i(x-x_0)$ with
$L_i\in\;\mbox{Hom}(\mathbb{R}^m,\mathbb{R}^n)$.
\end{remark}
\begin{definition} If $f$ is affinely approximatable at $x_0$, then\\
(a) $\sum_{i=1}^Q[[L_i]]\in
\mathbb{Q}_Q(\mbox{Hom}(\mathbb{R}^m,\mathbb{R}^n))$, denoted by
$Df(x_0)$ is defined as the differential of $f$ at $x_0$. We let
$|Df(x_0)|^2=\sum_{i=1}^Q|L_i|^2$, where $|L|$ is the Euclidean
norm of the
matrix associated with any $L\in\;\mbox{Hom}(\mathbb{R}^m,\mathbb{R}^n).$\\
(b) $\sum_{i=1}^Q[[L_i(v)]]$ is defined as the derivative of $f$
at $x_0$ in the direction $v$ and is denoted by $D_v f\in
\mathbb{Q}$. Let $|D_v f(x_0)|^2=\sum_{i=1}^Q |L_i(v)|^2$.
\end{definition}
In this section we will define derivative in a more natural way,
as we do in one variable calculus:
$$f'(x_0)=\lim_{x\to x_0}\DF{f(x)-f(x_0)}{x-x_0}.$$
The major difficulty is due to the fact that ``subtraction"
between two elements in $\mathbb{Q}$ is generally not well-defined
except for some special cases. For example, we can define the
subtraction $(-)$: $\mathbb{Q}\times \mathcal{I}_1\to \mathbb{Q}$,
for any $z=\sum_{i=1}^Q [[z_i]]\in \mathbb{Q}$ and $q=Q[[q_0]]\in
\mathcal{I}_1$ by setting $z(-)q=\sum_{i=1}^Q [[z_i-q_0]].$\\
Actually, as we will see, as long as two points are close enough
with each other, we can define ``subtraction'' between them by
matching corresponding components.
\begin{definition}
For $q=\sum_{i=1}^J k_i[[q_i]]\in \mathbb{Q}\sim \mathcal{I}_1$,
with $S(q)=(J,k_1,\cdot\cdot\cdot,k_J)$, define
\begin{equation*}
\begin{split}
\mathbb{P}(q;r)&=\{\sum_{j=1}^Q[[z_j]]:z_1,\cdot\cdot\cdot,z_Q\in
\mathbb{R}^n \;\mbox{with card}\{j:z_j\in
\mathbb{B}_r^n(q_i)\}=k_i\\
&\mbox{for
each}\;i=1,\cdot\cdot\cdot,J\;\mbox{and}\;r<r_0=2^{-1}\inf\{|q_i-q_j|:1\le
i<j\le J\}\}. \end{split} \end{equation*}\end{definition}
\begin{remark} Any element $z\in \mathbb{P}(q;r)$ can be expressed as
$\sum_{i=1}^J \sum_{j=1}^{k_i} [[z_j^{(i)}]],$ where $z_j^{(i)}\in
\mathbb{B}_r^n(q_i)$ for all fixed $i$ and
$j=1,2,\cdot\cdot\cdot,k_i$. \end{remark}  Now we are able to
define subtraction between $q\in \mathbb{Q}\sim \mathcal{I}_1$ and
$z\in \mathbb{P}(q;r)$ as follows:
\begin{definition} For $q\in \mathbb{Q}\sim \mathcal{I}_1$
and $z\in \mathbb{P}(q;r)$, define
$$q(-)z=\sum_{i=1}^J \sum_{j=1}^{k_i} [[q_i-z_j^{(i)}]],$$
and
$$z(-)p=\sum_{i=1}^J \sum_{j=1}^{k_i} [[z_j^{(i)}-q_i]].$$
\end{definition}
We are ready to introduce the quotient-based-derivative of
multiple-valued functions. Since if a multiple-valued function is
differentiable (in the sense of $\cite{af}$) at some point, then
it must be continuous there, we may just assume that the
multiple-valued function $f$ is continuous throughout its
domain.\\
For simplicity, let's first consider the one-dimensional domain.
We fix $x_0\in \mathbb{R}$. \\Case 1: $f(x_0)\in \mathcal{I}_1$.
\begin{definition} Suppose $f(x_0)\in \mathcal{I}_1$. $f$ is
said to be differentiable at $x_0$ if the following limit exists
(i.e, each component has a finite limit):
$$\lim_{x\to x_0}\DF{f(x)(-)f(x_0)}{x-x_0}.$$
If so, the limit is denoted as $f'(x_0)\in\mathbb{Q}$, called the
derivative of $f$ at $x_0$.
\end{definition}
\begin{remark}
It is easy to check that this definition is equivalent to the one
in Definition 5.1 (b).
\end{remark}
Case 2: $f(x_0)=\sum_{i=1}^J k_i[[q_i]]\in \mathbb{Q}\sim
\mathcal{I}_1$, with $S(f(x_0))=(J,k_1,\cdot\cdot\cdot,k_J)$. \\
Since $f$ is continuous at $x_0$, $f(x)\in \mathbb{P}(f(x_0);r)$
when $x$ is close enough with $x_0$. So the quotient
$$\DF{f(x)(-)f(x_0)}{x-x_0}$$
makes sense when $|x-x_0|$ small enough.
\begin{definition} Suppose $f(x_0)\in \mathbb{Q}\sim
\mathcal{I}_1$. $f$ is said to be differentiable at $x_0$ if the
following limit exists (i.e, each component has a finite limit):
$$\lim_{x\to x_0}\DF{f(x)(-)f(x_0)}{x-x_0}.$$
If so, the limit is denoted as $f'(x_0)\in\mathbb{Q}$, called the
derivative of $f$ at $x_0$.
\end{definition}
\begin{remark}
It is easy to check that this definition is also equivalent to the
one in Definition 5.1 (b).
\end{remark}
For higher dimensional domain, we define the directional
derivative as follows: \begin{definition} Suppose
$f:\mathbb{R}^m\to\mathbb{Q}$ is continuous, fix
$x_0\in\mathbb{R}^m$ and $v\in\mathbb{R}^m$, the directional
derivative of $f$ at $x_0$ in the direction $v$ is the following
limit if it exists:
$$L(v):=\lim_{t\to 0}\DF{f(x_0+tv)(-)f(x_0)}{t}.$$
\end{definition}
\subsection{A Regular Selection Theorem}
\begin{theorem}[$\cite{gj}$ Proposition 5.2]
Let $f:[a,b]\subset\mathbb{R}\to\mathbb{Q}$ be a continuous
function, then there exist continuous functions
$f_1,f_2,\cdot\cdot\cdot,f_Q:[a,b]\to\mathbb{R}^n$ such that
$f=\sum_{i=1}^Q [[f_i]]$ on $[a,b]$.
\end{theorem}
\begin{remark} Generally a continuous multiple-valued function
does not necessarily admit a continuous decomposition as noticed
by $\cite{gj}$. Another interesting question is whether the
decomposition shown above preserves the differentiability of $f$.
The following example suggests that being ``affinely
approximatable'' is not enough:
$$f:[-1,1]\to\mathbb{Q}_2(\mathbb{R}),f(x)=[[x]]+[[-x]].$$
It is affinely approximatable but not strongly affinely
approximatable at $0$. The decompositions $f_1=|x|,f_2=-|x|$
certainly are not differentiable at $0$. With the assumption of
being strongly affinely approximatable, we can prove:
\end{remark}
\begin{theorem}
If $f:[a,b]\subset\mathbb{R}\to\mathbb{Q}$ is a continuous
function and $x_0\in(a,b)$. Suppose $f$ is strongly affinely
approximatable at $x_0$, then there exist continuous functions
$f_1,f_2,\cdot\cdot\cdot,f_Q:[a,b]\to\mathbb{R}^n$ such that
$f=\sum_{i=1}^Q [[f_i]]$, and each $f_i$ is differentiable at
$x_0$.
\end{theorem}
\begin{proof}
Let $f_i$ be a decomposition from Theorem 5.1. We also let
$g_1,g_2,\cdot\cdot\cdot,g_Q\in \mathbb{A}(1,n)$ such that
$$f'(x_0)=Af(x_0)=\sum_{i=1}^Q [[g_i]]=g.$$
Consider $f(x_0)$ which is the same as both $\sum_{i=1}^Q
[[f_i(x_0)]]$ and $\sum_{i=1}^Q [[g_i(x_0)]]$. Suppose
$S(f(x_0))=(J,k_1,\cdot\cdot\cdot,k_J)$, and reordering $f_i's$ if
necessary such that
$$f(x)=\sum_{i=1}^Q [[f_i(x)]],$$
$$f_1(x_0)=f_2(x_0)=\cdot\cdot\cdot=f_{k_1}(x_0),$$
$$f_{k_1+1}(x_0)=f_{k_1+2}(x_0)=\cdot\cdot\cdot=f_{k_1+k_2}(x_0),$$
$$\cdot\cdot\cdot$$
$$f_{k_1+k_2+\cdot\cdot\cdot+k_{J-1}+1}(x_0)=f_{k_1+k_2+\cdot\cdot\cdot+k_{J-1}+2}(x_0)=\cdot\cdot\cdot=f_Q(x_0).$$
By the assumption
that $f$ is strongly affinely approximatable at $x_0$, we may
rewrite $g$ as
$$g=\sum_{i=1}^J k_i[[g_i]],$$
for $J$ non-crossing affine maps $g_i$.\\ By continuity of $f$,
apparently we have
\begin{equation*}
\begin{split}
\mathcal{G}^2(f(x),g(x))&=\sum_{i=1}^J
\mathcal{G}^2(\sum_{j=k_1+\cdot\cdot\cdot+k_{i-1}+1}^{k_1+\cdot\cdot\cdot+k_i}
[[f_j(x)]],k_i[[g_i(x)]])\\
&=\sum_{i=1}^J
\sum_{j=k_1+\cdot\cdot\cdot+k_{i-1}+1}^{k_1+\cdot\cdot\cdot+k_i}
|f_j(x)-g_i(x)|^2,
\end{split}
\end{equation*}
when $|x-x_0|$ small enough. The assumption that
$$\lim_{x\to x_0}|x-a|^{-1}\mathcal{G}(f(x),g(x))=0$$
shows that each $f_j$ is differentiable at $x_0$ and
$$f'_j(x_0)=g_i(x_0),j=k_1+\cdot\cdot\cdot+k_{i-1}+1,\cdot\cdot\cdot,k_1+\cdot\cdot\cdot+k_i,i=1,2,\cdot\cdot\cdot,J.$$
\end{proof}

\end{document}